\documentclass{article}[16pt]
\usepackage{amsfonts,amsmath}
\usepackage[russian]{babel}
\usepackage[cp1251]{inputenc}

\newtheorem{theorem}{Теорема}

\begin{document}

\title{О сходимости произведений \\
направленностей операторов}

\author{Д.В. Фуфаев \thanks{Механико-математический факультет Московского Государственного Университета имени М.В. Ломоносова,
Россия; e-mail:
fufaevdv@rambler.ru}  }

\date{}

\maketitle

Классические задачи анализа о сходимости тех или иных последовательностей зачастую сводятся к вопросу сходимости результатов действия некоторых операторов --- классическим примером является задача сходимости к исходной функции рядов Фурье или разнообразных средних, с рядами Фурье связанных; соответствующими операторами здесь являются операторы свертки с определенными функциями, например, с ядром Дирихле в случае частичных сумм ряда Фурье. В частности, в рамках данной работы нас будут интересовать результаты, связанные со средними Фейера и Абеля-Пуассона. Другим классическим примером, тесно связанным с первым, является задача дифференцирования неопределенного интеграла суммируемой функции; стоит отметить, что в этой задаче весьма плодотворным оказалось исследование максимального оператора, построенного по данному семейству операторов (см. [2, гл.1, \S1]). Нами рассматривается вопрос сходимости результатов действия упомянутых операторов к исходной функции почти всюду. В одномерном случае результат положителен для любой суммируемой функции. В случае же функций нескольких переменных ситуация усложняется: результат зависит от поведения индексов по разным координатам друг относительно друга. Если это поведение регулярно (см. [1, гл.XVII, \S 3]), то, как и в одномерном случае, достаточно суммируемости, причем факт регулярности может быть переформулирован в терминах соответствующего максимального оператора. В нерегулярной же ситуации приходится требовать принадлежность функции классу $L (\ln^+L)^{N-1}$, где $N$ --- размерность пространства (см. [1, гл.XVII, \S 2]).

Оказывается, что соответствующие результаты о сходимости почти всюду не используют специфики классического контекста. В данной работе приводится их обобщение для сходимости почти всюду в пространстве с абстрактной мерой. При этом в классической ситуации обнаруживаются новые возможности для определения регулярности, точнее - промежуточной регулярности, которые ранее не рассматривались.

Пусть $(X,\mu), (Y,\nu)$  ---  пространства с мерами, а оператор $T$ действует из $L^0(X,\mu)$ в $L^0(Y,\nu)$, где $L^0$ --- множество классов эквивалентности (эквивалентность --- совпадение почти всюду) измеримых действительнозначных (либо комплекснозначных) функций. Назовем $T$ сублинейным (или выпуклым), если из существования $Tf_1$ и $Tf_2$ следует существование $T(f_1+f_2)$ и при этом выполнено $|T(f_1+f_2)(x)|\le |Tf_1(x)|+|Tf_2(x)|$. Будем говорить, что выпуклый оператор $T$ имеет слабый тип $(1,1)$ с константой $C$, если для любого $\lambda>0$ и для любой функции $f\in L^1(X,\mu)$ выполняется следующее неравенство:

$$
\nu\{y\in Y: |Tf(y)|>\lambda\}\le \frac{C}{\lambda} ||f||_{L^1(X,\mu)}.
$$

В дальнейшем $(Y,\nu)=(X,\mu)$.


Для произвольной измеримой функции $g$ на $X$ определим функцию распределения \linebreak $\lambda_g(a)=\mu\{x\in X: |g(x)|>a \}$

Для нее нетрудно установить следующее равенство, аналогичное [3, гл.IX, $\S$4, следствие 2]:

\begin{equation}\label{Fuf_eq1}
\int\limits_{x:g(x)>\varepsilon} |g(x)|d\mu(x)=-\int\limits_{\varepsilon}^\infty a \  d\lambda_g(a).
\end{equation}

\textbf{Лемма 1.} \textit{Пусть $\mu(X)=1$, $T$ --- оператор слабого типа (1,1) с константой $C$, а функция $f$ такова, что $0\le f \le a$. Тогда для любого $\gamma>0$ найдутся число $k\in \mathbb N$, не зависящее от $f$, и множество $X_\gamma$ такое, что $\mu(X_\gamma)>1-\gamma$ и для любого $x\in X_{\gamma}$ выполнено
$$
 |Tf(x)|\le a\cdot 2^k
$$
}

 {\bfДоказательство.} Для любого натурального $n$ имеем
$$
\mu\{x\in X:|Tf(x)|>a\cdot 2^n\}\le \frac{C}{a\cdot 2^n}||f||_{L^1}\le \frac{C}{2^n}
$$
Для заданного $\gamma>0$ найдем номер $k$ такой, что $\frac{C}{2^k}<\gamma$. Тогда справедливо неравенство
$$
\mu\{x\in X:|Tf(x)|>a\cdot 2^k\}\le\gamma
$$
Следовательно, на множестве меры, не меньшей $1-\gamma$ выполняется неравенство $|Tf(x)|\le a\cdot 2^k$.

Лемма доказана

\textbf{Лемма 2.} \textit{ Пусть $\mu(X)=1$, $T$ --- оператор слабого типа (1,1) с константой $C$. Тогда для любого $\gamma>0$ найдутся число $k\in \mathbb N$ и множество $X_\gamma$ такое, что $\mu(X_\gamma)>1-\gamma$ и для любой функции $f$ и любого $a>0$ выполнено
$$
 \lambda_{Tf|_{X_\gamma}}(a)\le \frac{2\cdot C}{a } \int\limits_{x\in X_\gamma: |f(x)|>a\cdot 2^{-k-1}}|f(x)|d\mu(x)
$$
}

\textbf {Доказательство.} Для произвольного $b>0$ разделим $f$ на большую и малую части:

\[
f_1(x)=\begin{cases}
f(x),&\text{при $|f(x)|<b$;}\\
0,&\text{при $|f(x)|\ge b$,}
\end{cases}
\]
и $f_2(x)=f(x)-f_1(x).$

Очевидно, что

$$
|f(x)|\le |f_2(x)|+|f_1(x)|,
$$
из выпуклости оператора $T$ следует, что
$$
|Tf(x)|\le |Tf_2(x)|+|Tf_1(x)|,
$$
а из леммы 1 --- что существуют $k\in \mathbb N$ и $X_\gamma$ с $\mu(X_\gamma)>1-\gamma$ такие, что при $x\in X_\gamma$ выполнено
$$
|Tf(x)|\le |Tf_2(x)|+b\cdot 2^k,
$$
следовательно, справедливо следующее включение
$$
\{x\in X_\gamma: |Tf(x)|>b\cdot 2^{k+1}\}\subset \{x\in X_\gamma: |Tf_2(x)|>b\cdot 2^k\}
$$
откуда
$$
\lambda_{Tf|_{X_\gamma}}(b\cdot 2^{k+1})\le \lambda_{Tf_2|_{X_\gamma}}(b\cdot 2^{k})\le \frac{C}{b\cdot 2^k} ||f_2||_{L^1(X_\gamma,\mu)}=
\frac{C}{b\cdot 2^k} \int\limits_{x\in X_\gamma: |f(x)|>b}|f(x)|d\mu(x)
$$
положив $b=a\cdot 2^{-k-1}$, получаем
$$
\lambda_{Tf|_{X_\gamma}}(a)\le\frac{2\cdot C}{a} \int\limits_{x\in X_\gamma: |f(x)|>a\cdot 2^{-k-1}}|f(x)|d\mu(x).
$$
Лемма доказана.

\textbf{Теорема 1.}  (неравенство Харди-Литтльвуда). 

{\itПусть $\mu(X)=1$, $f\in L^1(X), f\cdot \ln(f+1)\in L^1(X), f \ge 0$, $T$ --- оператор слабого типа (1,1) и задано $\varepsilon>0$. Тогда для любого $\gamma>0$ найдется  множество $X_\gamma$ такое, что $\mu(X_\gamma)>1-\gamma$ и справедливо неравенство

$$
\int\limits_{X_\gamma}|Tf(x)|d\mu(x) \le A \int\limits_{X_\gamma}f(x)\cdot \ln(f(x)+1)d\mu(x) + B \int\limits_{X_\gamma}f(x)d\mu(x) + \varepsilon,
$$

где $A$ и $B$ --- постоянные, зависящие от $\varepsilon$, но не от $f$.}

\textbf{Доказательство.}

Множество $X_\gamma$ (а вдобавок --- и число $k\in \mathbb N$) дает лемма 2.

Оценим интеграл, разбивая функцию $Tf$ на ее большую и малую части и используя равенство \eqref{Fuf_eq1}:

$$
\int\limits_{X_\gamma}|Tf(x)|d\mu(x)= \int\limits_{x\in X_\gamma: |Tf(x)|\le\varepsilon }|Tf(x)|d\mu(x) + \int\limits_{x\in X_\gamma: |Tf(x)|>\varepsilon }|Tf(x)|d\mu(x)\le
$$

\begin{equation}\label{Fuf_eq15}
\le \varepsilon - \int\limits_{\varepsilon}^\infty a \  d\lambda_{Tf|_{X_\gamma}}(a)\le\varepsilon + \int\limits_{\varepsilon}^\infty \lambda_{Tf|_{X_\gamma}}(a)\  da +\lambda_{Tf|_{X_\gamma}}(\varepsilon)\cdot\varepsilon.
\end{equation}
В последнем неравенстве проинтегрировали по частям и использовали тот факт, что $\lambda \ge0$. Оценим слагаемые, используя лемму 2:

$$
\lambda_{Tf|_{X_\gamma}}(\varepsilon)\cdot\varepsilon \le  \frac {2\cdot C}{\varepsilon} \cdot \varepsilon \int\limits_{x: f(x)>\varepsilon\cdot 2^{-k-1}}f(x)d\mu(x) \le 2\cdot C\cdot ||f||_{L^1},
$$

$$
 \int\limits_{\varepsilon}^\infty \lambda_{Tf|_{X_\gamma}}(a)\  da \le  2\cdot C  \int\limits_{\varepsilon}^\infty  \frac{1}{a}\int\limits_{x: f(x)>a\cdot 2^{-k-1}}f(x)d\mu(x) da=
 $$
 
 $$
=2\cdot C\int\limits_{x: |f(x)|>\varepsilon\cdot 2^{-k-1}}  \int\limits_{\varepsilon}^{2^{k+1}f(x)} \frac{f(x)}{a}da d\mu(x)
 =2C\int\limits_{x: |f(x)|>\varepsilon\cdot 2^{-k-1}} f(x) ( \ln(2^{k+1}f(x)) - \ln\varepsilon     )      d\mu(x)\le 
 $$
 
 $$
\le 2C(||f\cdot \ln (f+1)||_{L^1}+((k+1)\ln2 +|\ln\varepsilon|)||f||_{L^1});
$$
подставляя эти оценки в  \eqref{Fuf_eq15}, получили нужное неравенство с постоянными $A=2\cdot C$ и $B=2\cdot C(1+(k+1)\cdot \ln2 +|\ln\varepsilon|)$.

\textbf {Замечание.} В частных случаях (максимальная функция\linebreak Харди-Литтлвуда, мажоранты средних Фейера) как правило оказывается, что неравенства, доказанные в леммах 1 и 2 и в теореме 1 следуют непосредственно из общего вида операторов, причем с $k=\gamma=0$. 

Порой в гармоническом анализе возникают приближения не последовательностью операторов, а семейством операторов, которые образуют лишь частично упорядоченное множество (т.е. не каждые два оператора сравнимы по номеру). Самый распространенный пример --- кратные ряды Фурье. Чтобы работать с такими семействами, вспомним понятие направленности (см. [4, т.2, стр.12]).

\textbf{Определение 1}. Непустое множество $A$ называется направленным, если на нем задан частичный порядок, удовлетворяющий следующему условию: для любых $m,n\in A$ найдется элемент $k\in A$ такой, что $m\le k$ и $n\le k$. Направленностью в множестве $X$ называется набор элементов $\{x_n\}_{n\in A}$, индексируемых элементами направленного множества. Направленность $\{x_n\}_{n\in A}$ в топологическом пространстве $X$ сходится к элементу $x$, если для любого непустого открытого множества $U$, содержащего $x$, найдется такой элемент $n_0\in A$, что $x_n\in U$ для всех $n\ge n_0, n\in A$. Понятным образом определяется сходимость числовых направленностей, а также поточечная сходимость и сходимость почти всюду направленностей числовых функций.

Пусть $\{T_n\}_{n\in A}$ --- направленность линейных операторов, переводящих $L^0(X,\mu)$ в себя. Максимальным оператором относительного данного семейства операторов будем называть оператор $T:f(x)\mapsto \sup\limits_{n\in A}|T_nf(x)|$. Будем рассматривать лишь счетные направленности, т.е. такие, что множество $A$ счетно --- это гарантирует измеримость функции $Tf(x)$. Нетрудно увидеть, что оператор $T$ является сублинейным.

Аналогично [5, теорема 5.1.3] доказывается следующий результат.

\textbf{Теорема 2.}\textit{  Пусть направленность линейных операторов $\{T_n\}_{n\in A}$ такова, что соответствующий максимальный оператор $T$ имеет слабый тип (1,1), и пусть для любой функции $\phi$ из всюду плотного в $L^1(X,\mu)$ множества $\lim\limits_{n\in A}T_n\phi(x)=\phi(x)$ почти всюду на $X$. Тогда для любой функции $f\in L^1(X,\mu)$ выполняется $\lim\limits_{n\in A}T_nf(x)=f(x)$ почти всюду на $X$.
}

В дальнейшем нас будут интересовать направленности лишь интегральных операторов, то есть операторов вида $Tf(x)=\int\limits_X K(x,y)f(y)d\mu(y)$.

Теперь нам нужно определить функции и операторы на них для произведений пространств. Для работы с пространством функций вспомним понятие тензорного произведения [6, 10.42.1, 10.42.2]: пусть $f^1\in L^1(X^1,\mu^1), f^2\in L^1(X^2,\mu^2)$, тогда их тензорным произведением называется функция $f^1\otimes f^2\in L^1(X^1\times X^2, \mu^1\otimes \mu^2)$, определяемая по формуле $f^1\otimes f^2(x_1,x_2)=f^1(x_1)f^2(x_2)$, где $\mu^1\otimes \mu^2$ --- это произведение мер $\mu^1$ и $\mu^2$ (см., например, [4, т.1, стр. 223]). Такие функции называются элементарными тензорами. Алгебраическим тензорным произведением $L^1(X^1,\mu^1)$ и $L^1(X^2,\mu^2)$ называется пространство линейных комбинаций элементарных тензоров и обозначается $L^1(X^1,\mu^1) \mathbin{\otimes} L^1( X^2,\mu^2)$. Проективным тензорным произведением, обозначаемым $L^1(X^1,\mu^1) \mathbin{\hat{\otimes}} L^1( X^2,\mu^2)$, называется пополнение алгебраического по 
так называемой проективной норме [7, 0.3.28], но нам она не нужна, а нужен тот факт, что оно
изометрически изоморфно $L^1(X^1\times X^2,\mu^1\otimes \mu^2)$ [7, следствие 0.3.36]. В частности, важен тот факт, что произвольная функция из $L^1(X^1\times X^2,\mu^1\otimes \mu^2)$ приближается комбинациями элементарных тензоров, т.е. функциями вида $\sum\limits_{k=1}^n f_k(x_1)\cdot g_k(x_2)$. Более того, понятно, что в качестве $f_k$ и $g_k$ можно брать функции из всюду плотного в соответствующем $L^1(X^i)$ множества, например, из подпространства, состоящего из ограниченных измеримых функций, которое мы будем обозначать $L^{\infty}(X^i)$, что, конечно, не вызовет путаницы.

Вспомним, что функция $\phi(x)=x\ln^{D}(x+1)$ - выпукла для $x\ge0$ и $D\in \mathbb N$.


Для измеримого пространства $(X,\mu)$ через $L(\ln^{+}L)^D  (X), D\in \mathbb N,$ обозначается пространство измеримых комплекснозначных функций $f$, удовлетворяющих условию

$$
\int\limits_{X}|f(x)|\ln^D(|f(x)|+1)d\mu(x)<\infty.
$$

\textbf{Лемма 3.} \textit{ Пусть $(X^1,\mu^1)$, $(X^2,\mu^2)$ --- пространства конечной меры. Тогда $L^{\infty}(X^1)\otimes L^{\infty}(X^2)$ всюду плотно в $L\ln^{+} L (X^1\times X^2)$ относительно величины $\rho(f,g)=\int\limits_{X}|f(x)-g(x)|\ln(|f(x)-g(x)|+1)d\mu(x)$, где $\mu=\mu^1\otimes\mu^2$.
}

\textbf{Доказательство.} Возьмем произвольную функцию $f\in L\ln^{+} L (X^1\times X^2)$ и зафиксируем произвольное $\varepsilon>0$. В силу неравенства Чебышева выполнено $\mu\{x\in X:|f(x)|\ge n\}\le\frac{1}{n}||f||_{L^1}$. Обозначим $B_{n,k}=\{\frac{k}{n}\le f(x)<\frac{k+1}{n}\}$ и $B_n=\bigsqcup\limits_{k=-n^2}^{n^2-1}B_{n,k}\subset\{x:|f(x)|<n+1\}$. Тогда в силу абсолютной непрерывности интеграла Лебега найдется такое натуральное $n$, что $\int\limits_{X\setminus B_n}|f(x)|\ln(|f(x)|+1)d\mu(x)<\varepsilon$. Пусть, к тому же, $n$ удовлетворяет условию $\frac{\mu(X)}{n}<\sqrt{\varepsilon}$. Далее, по определению меры на произведении, для произвольного $\delta>0$ найдутся такие множества $D_{n,k,i}, i=1\dots, r_k$, имеющие вид произведений измеримых множеств из $X^1$ и $X^2$, которые можно взять дизъюнктными при фиксированных $n$ и $k$, что $\mu(B_{n,k}    \bigtriangleup \bigsqcup\limits_{i=1}^{r_k}D_{n,k,i})<\frac{\delta}{2n^6}$. Положим $D_n=\bigcup\limits_{k,i}D_{n,k,i}$ и заметим, что $\mu(B_n \bigtriangleup D_n)\le 2 n^2\frac{\delta}{2n^6}=\frac{\delta}{n^4}$. Пользуясь теперь абсолютной непрерывностью, возьмем $\delta$ соответствующим $\varepsilon$ для функций $(|f(x)|+n)\ln(|f(x)|+n+1)$ и $|f(x)|\ln(|f(x)|+1)$, и вдобавок $\delta<\varepsilon$. Пусть $\varphi_n(x)=\sum\limits_{k=-n^2}^{n^2-1}\sum\limits_{i=1}^{r_k}\frac{k}{n}\chi_{D_{n,k,i}}(x)$, $h_n(x)=\sum\limits_{k=-n^2}^{n^2-1}\frac{k}{n}\chi_{B_{n,k}}(x)$ и оценим $\rho(f,\varphi_n)$:

$$
\int\limits_{X}|f(x)-\varphi_n(x)|\ln(|f(x)-\varphi_n(x)|+1)d\mu(x)
=\int\limits_{X\setminus D_n}|f(x)-\varphi_n(x)|\ln(|f(x)-\varphi_n(x)|+1)d\mu(x)+
$$

$$
+\int\limits_{D_n\cap B_n}|f(x)-\varphi_n(x)|\ln(|f(x)-\varphi_n(x)|+1)d\mu(x)+
\int\limits_{D_n\setminus B_n}|f(x)-\varphi_n(x)|\ln(|f(x)-\varphi_n(x)|+1)d\mu(x)\le
$$

$$
\le\int\limits_{X\setminus D_n}|f(x)|\ln(|f(x)|+1)d\mu(x)+\int\limits_{D_n\cap B_n}|f(x)-\varphi_n(x)|^2d\mu(x)+
$$

$$
+\int\limits_{D_n\setminus B_n}(|f(x)|+n)\ln(|f(x)|+n+1)d\mu(x)\le
$$

$$
\le\int\limits_{(X\setminus B_n)\cup(B_n\bigtriangleup D_n)}|f(x)|\ln(|f(x)|+1)d\mu(x)+||f(x)-\varphi_n(x)||^2_{L^2(D_n\cap B_n)}+
$$

$$
+\int\limits_{B_n\bigtriangleup D_n}(|f(x)|+n)\ln(|f(x)|+n+1)d\mu(x)\le
$$

$$
\le2\varepsilon+     \left(||f(x)-h_n(x)||_{L^2(D_n\cap B_n)}+||h_n(x)-\varphi_n(x)||_{L^2(D_n\cap B_n)}\right)^2  +\varepsilon\le
$$

$$
\le3\varepsilon+ \Biggl(\sqrt{\int\limits_{D_n\cap B_n} (f(x)-\sum\limits_{k=-n^2}^{n^2-1}\frac{k}{n}\chi_{B_{n,k}}(x))^2d\mu(x)}
+
||\sum\limits_{k=-n^2}^{n^2-1}\frac{k}{n}(\chi_{B_{n,k}}-\sum\limits_{i=1}^{r_k}\chi_{D_{n,k,i}})||_{L^2(D_n\cap B_n)}\Biggr)^2=
$$

$$
=3\varepsilon+ \Biggl(\sqrt{\int\limits_{D_n\cap B_n} (\frac{1}{n})^2d\mu(x)}+\sum\limits_{k=-n^2}^{n^2-1}\frac{|k|}{n}||(\chi_{B_{n,k}}(x)-
\sum\limits_{i=1}^{r_k}\chi_{D_{n,k,i}}(x))||_{L^2(D_n\cap B_n)}\Biggr)^2\le
$$

$$
\le3\varepsilon+\left(\frac{\sqrt{\mu(X)}}{n}+\sum\limits_{k=-n^2}^{n^2-1}\frac{n^2}{n}\sqrt{\int\limits_{D_n\cap B_n}(\chi_{B_{n,k}}(x)-\sum\limits_{i=1}^{r_k}\chi_{D_{n,k,i}}(x))^2d\mu(x)}\right)^2\le
$$



$$
\le3\varepsilon+\left(\sqrt{\varepsilon}+\sum\limits_{k=-n^2}^{n^2-1}n\sqrt{\mu(B_{n,k}\bigtriangleup\bigsqcup\limits_{i=1}^{r_k}D_{n,k,i}) }\right)^2\le3\varepsilon+\left(\sqrt{\varepsilon}+2n^3\sqrt{\frac{\delta}{2n^6}) }\right)^2
\le3\varepsilon+\varepsilon\left(1+\frac{2}{\sqrt{2}}) \right)^2\le12\varepsilon, 
$$
откуда, учитывая, что $\varphi_n$ имеет вид линейной комбинации элементарных тензоров, следует требуемое утверждение.

Лемма доказана.







Далее, если линейные операторы $T^i$ заданы на $L^1(X^i,\mu^i)$, $i=1,2$, то существует единственный оператор $T^1\mathbin{\hat{\otimes}}T^2$, действующий в $L^1(X^1,\mu^1) \mathbin{\hat{\otimes}} L^1( X^2,\mu^2)$, такой, что на элементарных тензорах он действует по формуле: $T^1\mathbin{\hat{\otimes}}T^2[f\cdot g](x,y)=T^1f(x)\cdot T^2g(y)$ [7, теорема 0.3.40].
Очевидно, что тензорным произведением двух интегральных операторов с ядрами $K_i(\cdot,\cdot)$ будет интегральный оператор с ядром $K_1(\cdot,\cdot)\cdot K_2(\cdot,\cdot)$.

Для двух направленностей операторов $\{T^1_{n^1}\}_{n^1\in A^1}$ и $\{T^2_{n^2}\}_{n^2\in A^2}$ их тензорным произведением назовем направленность
$\{T^1_{n^1} \mathbin{\hat{\otimes}}T^2_{n^2} \}_{\substack {   \textbf n\in A }} $, где $\textbf n=(n^1,n^2), A=A^1\times A^2$, причем $(n^1,n^2)>(m^1,m^2)$ тогда и только тогда, когда $n^1>m^1$ и $n^2>m^2$.
Тензорное произведение большего числа множителей определяется очевидным образом.

\textbf{Теорема 3.} {\it Пусть $(X^i,\mu^i)$, $i=1,\dots,D$ --- измеримые пространства меры 1 каждое, $X=\prod\limits_{i=1}^DX^i$, $\mu=\bigotimes\limits_{i=1}^D\mu^i$, $\{T^i_{n^i}\}_{n^i\in A^i}$, $i=1,\dots,D$ --- направленности линейных интегральных операторов (т.е., имеющих вид $T^i_{n^i}f(x)=\int\limits_{X^i} K^i_{n^i}(x,u) f(u,y)d\mu^{i}(u)$), действующих в соответстующих $L^1(X^i,\mu^i)$, таких, что каждый максимальный оператор $T^i$ имеет слабый тип (1,1) и, кроме того, для любого $i=1,\dots,D$ и любой ограниченной функции $\phi\in L^1(X^i,\mu^i)$ выполнено $\lim\limits_{n^i\in A^i}T^i_{n^i}\phi(x)=\phi(x)$ $\mu^i-$почти всюду. 

Тогда для любой $f\in L(\ln^{+})^D L(X)$ выполнено $\lim\limits_{\textbf n\in A}T_{\textbf n}f(x)=f(x)$ $\mu-$почти всюду, где $A$ --- произвольная поднаправленность тензорного произведения направленностей $\{T^i_{n^i}\}_{n^i\in A^i}$.
}

\textbf{Доказательство.} Докажем индукцией по $D$.  Случай $D=1$ --- результат теоремы 2. Пусть $D>1$. Предположим, что для всех $k<D$ утверждение доказано, докажем для $k=D$. 

Очевидно, можно предположить, что функция $f$ неотрицательна.

Обозначим точки из $X$ как $(x,y)$, где $x\in X^1, y\in Y=\prod\limits_{i=2}^DX^i$, тогда функцию $f$ можно обозначать как $f(x,y)$. Кроме того, через $\mu^{1\bot}$ будем обозначать меру $\bigotimes\limits_{i=2}^D\mu^i$, а через $n^{1\bot}$ --- элемент из $A^2\times\dots\times A^D$

Возьмем произвольное $\gamma>0$ и, пользуясь теоремой 1, найдем множества $X^i_\gamma$, $i=1,\dots,D$, такие, что $\mu^1(X^1_\gamma)>1-\gamma$ и $\mu^{1\bot}(Y_\gamma)=\mu^{1\bot}(\prod\limits_{i=2}^DX^i_\gamma)>(1-\gamma)^{D-1}>1-(D-1)\gamma$. 

Кроме того, по теореме Егорова (которая, как нетрудно видеть, остается в силе для счетных направленностей), в силу того, что, по условию $\lim\limits_{n^1\in A^1}    \int\limits_{X^1}K^1_{n^1}(x,t)d\mu^1(t)=1$ $\mu^1-$почти всюду, можно выбрать такое множество $X^1_{2\gamma}\subset X^1_\gamma$, что $\int\limits_{X^1}K^1_{n^1}(x,t)d\mu^1(t)$ сходится к 1 равномерно на $X^1_{2\gamma}$, при этом $\mu^1(X^1_{2\gamma})>1-2\gamma.$ Через $X_\gamma$ будем обозначать $X^1_{2\gamma}\times Y_\gamma.$

Далее рассматриваем ограничения всех функций и операторов на эти множества.

Положим $g_n(x,y)$ --- функции вида $\sum\limits_{k=1}^ng_k^1(x)\cdot g_k^{1\bot}(y)$, приближающие функцию $f$ по величине $\rho(\cdot,\cdot)$ на пространстве $X_\gamma$, введенной в лемме 3 (а значит, и по норме $||\cdot||_{X_\gamma}$). Далее, $h_n=f-g_n$. Возьмем произвольное $\varepsilon>0$. Построим функции $T^1h_n$ по первым координатам, а именно

$$
T^1h_n (x,y) =\sup_{n^1 \in A^1} T^1_{n^1} h_n(x,y)
$$

По Теореме 1 справедливо неравенство:

$$
\int\limits_{Y_\gamma} \int\limits_{X^1_{2\gamma}} T^1h_n(x,y)d\mu^1(x)d\mu^{1\bot}(y) \le 
$$

$$
\le A\int\limits_{Y_\gamma} \int\limits_{X^1_{2\gamma}}|h_n(x,y)|\cdot \ln(|h_n(x,y)|+1)d\mu^{1}(x)d\mu^{1\bot}(y) +
 B\int\limits_{Y_\gamma} \int\limits_{X^1_{2\gamma}}|h_n(x,y)|d\mu^{1}(x)d\mu^{1\bot}(y) + \frac{1}{2}\varepsilon^2,
$$
из которого, в частности, следует, что $T^1h_n(x,y)\in L^1 (X_\gamma)$.

Интегралы в правой части неравенства стремятся к нулю, поэтому существует такое число $K$, что (для краткости $h=h_K, g=g_K$)  $\int\limits_{Y_\gamma} \int\limits_{X^1_{2\gamma}} T^1h (x,y)d\mu^{1}(x)d\mu^{1\bot}(y)<\varepsilon^2$.

Пусть $E$ - множество таких точек $(x_0,y_0)\in X_\gamma$, что

1.  $\int\limits_{Y_\gamma} T^1h(x_0,t)\ln^{D-2}(T^1h(x_0,t)+1)d\mu^{1\bot}(t) < \infty  $,

2. Направленность $T^{1\bot}_{n^{1\bot}}(T^1h)(x_0,y_0)$ сходится к $T^1h(x_0,y_0)$

Ниже мы докажем, что  $\int\limits_{X_\gamma} T^1h(x,y) \ln^{D-2}(T^1h(x,y)+1) d\mu^1(x)d\mu^{1\bot}(y) < \infty$ (заметим, что в случае $D=2$ это утверждение просто следует из  неравенства Харди-Литлвуда, поэтому доказательство будет проводиться для $D>2$). Тогда множество $E$ будет иметь меру $\mu(X_\gamma)$, так как по теореме Фубини первое условие будет выполняться для п.в. $x_0\in X^1_\gamma$ , а условие 2 выполняется для п.в. $y_0\in Y_\gamma$, как только выполнено условие 1, по предположению индукции.



Вспомним следующее неравенство Йенсена для интеграла (см. [8, гл. 3, теорема 3.3]): для выпуклой функции $\phi$ и интегрируемой функции $f$ справедливо неравенство

$$
\phi \left(\int\limits_{X}f(x)\frac{d\nu(x)}{\nu(X)} \right) \le \int\limits_{X}\phi(f(x))\frac{d\nu(x)}{\nu(X)},
$$
из которого следует справедливость следующего

$$
\phi \left(\int\limits_{X}f(x)d\nu(x) \right) \le\frac{1}{\nu(X)} \int\limits_{X}\phi(\nu(X)\cdot f(x))d\nu(x).
$$


По определению максимального оператора $T^1$, для любого $\varepsilon >0$ существует элемент направленности $n^1\in A^1$ (зависящий, вообще говоря, от $(x,y)$) такой, что

$$
 \int\limits_{X_\gamma} T^1h(x,y) \ln^{D-2}(T^1h(x,y)+1) d\mu^{1}(x)d\mu^{1\bot}(y) \le
 $$
 
 $$\le
 \int\limits_{X_\gamma}  (T^1_{n^1} h(x,y)+\varepsilon) \ln^{D-2}(T^1_{n^1} h(x,y) +\varepsilon+1) d\mu^{1}(x)d\mu^{1\bot}(y) =
$$

$$
=\int\limits_{X_\gamma}  \left(\int\limits_{X^1} K^1_{n^1}(x,u) h(u,y)d\mu^{1}(u)+\varepsilon\right)
 \ln^{D-2}\left(\int\limits_{X^1}K^1_{n^1}(x,u)h(u,y)d\mu^{1}(u) +\varepsilon+1\right) d\mu^{1}(x)d\mu^{1\bot}(y) =
$$

$$
=\int\limits_{X_\gamma}  \int\limits_{X^1}  K^1_{n^1}(x,u)\biggl( h(u,y)+\frac{\varepsilon}{\int\limits_{X^1}K^1_{n^1}(x,t)d\mu^1(t)} \biggl)d\mu^{1}(u)\times
$$
 
 $$
 \times\ln^{D-2}\left(\int\limits_{X^1}  K^1_{n^1}(x,u)\biggl( h(u,y)+\frac{\varepsilon}{\int\limits_{X^1}K^1_{n^1}(x,t)d\mu^1(t)} \biggl)d\mu^{1}(u) +1\right) d\mu^{1}(x)d\mu^{1\bot}(y) \le
$$
применим неравенство Йенсена для $\phi(x)=x\cdot\ln^{D-2}(x+1)$,  $f(u)= h(u,y)+\frac{\varepsilon}{\int\limits_{X^1}K^1_{n^1}(x,t)d\mu^1(t)}$, очевидно, интегрируемой, и $\nu(X^1)=\int\limits_{X^1}K^1_{n^1}(x,u)d\mu^1(u)$:

$$
\le\int\limits_{X_\gamma}    \frac{1}{\int\limits_{X^1}K^1_{n^1}(x,t)d\mu^1(t)}   \int\limits_{X^1}K^1_{n^1}(x,u)  \int\limits_{X^1}K^1_{n^1}(x,t)d\mu^1(t)  \biggl( h(u,y)+\frac{\varepsilon}{\int\limits_{X^1}K^1_{n^1}(x,t)d\mu^1(t)}\biggl)\times
$$

$$
\times\ln^{D-2}\left(  \int\limits_{X^1}K^1_{n^1}(x,t)d\mu^1(t) \biggl(h(u,y)+\frac{\varepsilon}{\int\limits_{X^1}K^1_{n^1}(x,t)d\mu^1(t)}\biggl)+1\right)d\mu^1(u)d\mu^1(x)d\mu^{1\bot}(y)\le
$$

$$
\le\int\limits_{X_\gamma}      \int\limits_{X^1}K^1_{n^1}(x,u)   \biggl( h(u,y)+\frac{\varepsilon}{\int\limits_{X^1}K^1_{n^1}(x,t)d\mu^1(t)}\biggl)\times
$$

$$
\times\ln^{D-2}\left(  \int\limits_{X^1}K^1_{n^1}(x,t)d\mu^1(t) \cdot h(u,y)+\varepsilon+1\right)d\mu^1(u)d\mu^1(x)d\mu^{1\bot}(y)\le
$$
 в силу равномерной сходимости $\int\limits_{X^1}K^1_{n^1}(x,t)d\mu^1(t)$ к 1 на $X^1_{2\gamma}$, найдется такой элемент направленности $n^1$, что, вдобавок, для всех $x\in X^1$ справедливо неравенство $ 1/2 <\int\limits_{X^1}K^1_{n^1}(x,t)d\mu^1(t)<1+\varepsilon$




$$
\le\int\limits_{X_\gamma}      \int\limits_{X^1}K^1_{n^1}(x,u)   \left( h(u,y)+2\varepsilon\right)\times\ln^{D-2}\left(  (1+\varepsilon) \cdot h(u,y)+\varepsilon+1\right)d\mu^1(u)d\mu^1(x)d\mu^{1\bot}(y)=
$$

$$
=\int\limits_{X_\gamma}      \int\limits_{X^1}K^1_{n^1}(x,u)   \left( h(u,y)+2\varepsilon\right)\times\ln^{D-2}\left(  (1+\varepsilon) \cdot (h(u,y)+1\right))d\mu^1(u)d\mu^1(x)d\mu^{1\bot}(y)\le
$$

$$
\le\int\limits_{X_\gamma}   \int\limits_{X^1}K^1_{n^1}(x,u) \left( h(u,y)+2\varepsilon\right)\cdot 2^{D-3}\cdot (\ln^{D-2}(h(u,y)+1)+ \ln^{D-2} (1+\varepsilon)  )d\mu^1(u)d\mu^1(x)d\mu^{1\bot}(y)=
$$

$$
=2^{D-3}\cdot\int\limits_{X_\gamma}   \int\limits_{X^1}K^1_{n^1}(x,u) \left( h(u,y)+2\varepsilon\right)\cdot  \ln^{D-2}(h(u,y)+1)d\mu^1(u)d\mu^1(x)d\mu^{1\bot}(y)+
$$

$$
+ 2^{D-3}\cdot \ln^{D-2} (1+\varepsilon) \cdot\int\limits_{X_\gamma}   \int\limits_{X^1}K^1_{n^1}(x,u) \left( h(u,y)+2\varepsilon\right)d\mu^1(u)d\mu^1(x)d\mu^{1\bot}(y)\le
$$

$$
\le2^{D-3}\int\limits_{X_\gamma}   T^1 [\left( h(x,y)+2\varepsilon\right)  \ln^{D-2}(h(x,y)+1)]d\mu^1(x)d\mu^{1\bot}(y)+
$$
$$
+2^{D-3} \ln^{D-2} (1+\varepsilon) \int\limits_{X_\gamma}   T^1 \left[ h(x,y)+2\varepsilon\right]d\mu^1(x)d\mu^{1\bot}(y)
$$
по неравенству Харди-Литтлвуда, последний интеграл конечен. В силу произвольности $\varepsilon>0$ и неравенства Харди-Литтлвуда

$$
 \int\limits_{X_\gamma} T^1h(x,y) \ln^{D-2}(T^1h(x,y)+1) d\mu^{1}(x)d\mu^{1\bot}(y) \le 2^{D -3}\cdot\int\limits_{X_\gamma}   T^1 [\left( h(x,y)\right)\cdot  \ln^{D-2}(h(x,y)+1)]d\mu^1(x)d\mu^{1\bot}(y)\le
$$


\begin{multline*}
\le 2^{D-3}(A \int\limits_{X_\gamma} h(x,y)\cdot\ln^{D-2}( h(x,y)+1)\cdot\ln [h(x,y) \ln^{D-2}( h(x,y)+1)+1]    \ d\mu^1(x)d\mu^{1\bot}(y)+\\
+ B \int\limits_{X_\gamma} h(x,y)\cdot \ln^{D-2}( h(x,y)+1)\  d\mu^1(x)d\mu^{1\bot}(y)  + \gamma);
\end{multline*}
второй интеграл, очевидно, сходится. Оценим первый интеграл:

$$
\int\limits_{X_\gamma} h(x,y)\cdot \ln^{D-2}( h(x,y)+1)\cdot    \ln [h(x,y) \ln^{D-2}( h(x,y)+1)+1]    \ d\mu^1(x)d\mu^{1\bot}(y) \le
$$

$$
\le\int\limits_{X_\gamma} h(x,y)\cdot \ln^{D-2}( h(x,y)+1)   \ln [h(x,y) \ln^{D-2}( h(x,y)+1) +h(x,y)+  \ln^{D-2}( h(x,y)+1)   +1]    \ d\mu^1(x)d\mu^{1\bot}(y) =
$$

$$
=\int\limits_{X_\gamma} h(x,y)\cdot\ln^{D-2}( h(x,y)+1)\cdot    \ln [ (h(x,y)+1)   (\ln^{D-2}( h(x,y)+1) +1)]    \ d\mu^1(x)d\mu^{1\bot}(y)=
$$

$$
=\int\limits_{X_\gamma} h(x,y)\cdot\ln^{D-1}( h(x,y)+1)d\mu^1(x)d\mu^{1\bot}(y)+
$$

$$
+\int\limits_{X_\gamma} h(x,y)\cdot\ln^{D-2}( h(x,y)+1)\cdot  \ln [  \ln^{D-2}( h(x,y)+1) +1]    \ d\mu^1(x)d\mu^{1\bot}(y);
$$
первый из этих интегралов сходится по условию теоремы. Оценим второй интеграл:

$$
\int\limits_{X_\gamma} h(x,y)\cdot \ln^{D-2}( h(x,y)+1)\cdot \ln [  \ln^{D-2}( h(x,y)+1) +1]    \ d\mu^1(x)d\mu^{1\bot}(y)\le
$$

$$
\le\int\limits_{X_\gamma} h(x,y)\cdot \ln^{D-2}( h(x,y)+1)\cdot  \ln [  (\ln( h(x,y)+1) +1)^{D-2} ]    \ d\mu^1(x)d\mu^{1\bot}(y)=
$$

$$
=(D-2)\int\limits_{X_\gamma} h(x,y)\cdot \ln^{D-2}( h(x,y)+1)\cdot  \ln [  \ln( h(x,y)+1) +1 ]    \ d\mu^1(x)d\mu^{1\bot}(y)\le
$$

$$
\le(D-2)\int\limits_{X_\gamma} h(x,y)\cdot \ln^{D-2}( h(x,y)+1) \cdot \ln( h(x,y)+1)d\mu^1(x)d\mu^{1\bot}(y)=
$$

$$
=(D-2)\int\limits_{X_\gamma} h(x,y)\cdot \ln^{D-1}( h(x,y)+1)d\mu^1(x)d\mu^{1\bot}(y),
$$
а этот интеграл сходится. Таким образом, показали, что \newline $\int\limits_{X_\gamma} T^1h(x,y)\cdot \ln^{D-2}(T^1h(x,y)+1) d\mu^1(x)d\mu^{1\bot}(y) < \infty  $. Следовательно, мера множества $E$ равна $\mu(X_\gamma)$.


Пусть $(x_0,y_0)$ - точка множества $E$ и $T^1_{n^1} \mathbin{\hat{\otimes}}T^{1\bot}_{n^{1\bot}}$ --- произвольный оператор из направленности $A$. Имеем

$$
|T^1_{n^1} \mathbin{\hat{\otimes}}T^{1\bot}_{n^{1\bot}}h(x_0,y_0)|\le\int\limits_{Y}\int\limits_{X^1}|K^1_{n^1}(x_0,u)K^{1\bot}_{n^{1\bot}}(y_0,v)h(u,v)|d\mu^1(u)d\mu^{1\bot}(v)=
$$

$$
=\int\limits_{Y}|K^{1\bot}_{n^{1\bot}}(y_0,v)|\int\limits_{X^1}|K^1_{n^1}(x_0,u)h(u,v)|d\mu^1(u)d\mu^{1\bot}(v)\le  \int\limits_{Y}|K^{1\bot}_{n^{1\bot}}(y_0,v)|T^1h(x_0,v)d\mu^{1\bot}(v)
$$


Беря предел по направленности $A$, получаем $ {\varlimsup\limits_{\textbf n\in A}}T_{\textbf n}h(x_0,y_0) \le T^1h(x_0,y_0)$. Таким образом, так как $(x_0,y_0)$ - произвольная точка множества меры $\mu(X_\gamma)$, то получаем $0\le  {\varlimsup\limits_{\textbf n\in A}}T_{\textbf n}h(x_0,y_0) \le \varepsilon$ всюду в $X_\gamma$, кроме, быть может, множества меры меньшей, чем $\varepsilon$ (например, по неравенству Чебышева). С другой стороны, так как функция $g$ ограничена и имеет вид комбинации элементарных тензоров, имеем

$$
T^1_{n^1} \mathbin{\hat{\otimes}}T^{1\bot}_{n^{1\bot}}g(x_0,y_0)=T^1_{n^1} \mathbin{\hat{\otimes}}T^{1\bot}_{n^{1\bot}}\left( \sum\limits_{k=1}^Kg_k^1(x)\cdot g_k^{1\bot}(y)\right)=\sum\limits_{k=1}^KT^1_{n^1}   g_k^1(x)\cdot      T^{1\bot}_{n^{1\bot}} g_k^{1\bot}(y)\xrightarrow[\textbf n\in A]{}\sum\limits_{k=1}^Kg_k^1(x)\cdot g_k^{1\bot}(y)
$$
для почти всех $(x,y)\in X$. Таким образом, имеем


\begin{equation}\label{Fuf_eq2}
0\le  {\varlimsup\limits_{\textbf n\in A}}T_{\textbf n}f(x_0,y_0) -  {\varliminf\limits_{\textbf n\in A}}T_{\textbf n}f(x_0,y_0)  \le\varepsilon
\end{equation}
для точек из множества меры большей, чем $(1-\gamma)^{D-1}(1-2\gamma)-\varepsilon$. В силу произвольности $\varepsilon$ и $\gamma$, неравенство \eqref{Fuf_eq2} обращается в равенство почти всюду.

Теорема доказана.

\textbf {Замечание 1.} В теореме 3 рассматриваются не произвольные семейства линейных операторов, как в теореме 2, которую она обобщает, а лишь интегральные, по той причине, что в доказательстве для случая $D>2$ существенно используется  интегральное неравенство Йенсена. В случае же $D=2$ можно действительно брать произвольные линейные операторы, как и в теореме 2.

\textbf {Замечание 2.} Также в теореме 3 более сильное условие на всюду плотное множество, на котором изначально должна сходиться направленность операторов: произвольное всюду плотное подмножество $L^1$ может вообще не лежать в $L\ln^+L$,  в то время как в шаге индукции возникает необходимость приближения в терминах этого пространства --- то есть, по величине, определенной перед леммой 3. Однако, как видно из доказательства леммы 3, достаточно требовать сходимости не на всех ограниченных функциях, а только на ступенчатых. Более того, видно, что достаточно требовать сходимости на подмножествах, удовлетворяющих утверждению леммы 3 --- то есть, на таких подмножествах, произведения которых приближают функции из $L\ln^+ L$ на произведении, причем для каждого пространства-сомножителя может быть свое подмножество функций, на котором требуется сходимость в условии. Например, если какое-то из $X^i$ --- нормальное топологическое пространство (т.е. такое, в котором замкнутые множества отделяются открытыми окрестностями) с регулярной борелевской мерой (т.е. измеримые множества приближаются замкнутыми вписанными и открытыми описанными, таковыми например являются радоновские меры и в частности мера Лебега в $\mathbb R^n$), то, пользуясь большой леммой Урысона, по аналогии с [9, теорема 7.1.2] можно требовать сходимости на множестве непрерывных (или даже ограниченных непрерывных) на $X^i$ функций. Пользуясь же теоремами типа Уитни о продолжении для гладких функций (см., например, [10, теорема 3.1.14]), можно ограничиваться и гладкими функциями.

\textbf {Замечание 3.} Очевидно, изложенные результаты остаются верны для пространств конечной меры.

Простейшее применение полученного результата относится к теории суммирования кратных рядов Фурье. А именно, пусть $M_i$ --- натуральные числа ($i=1,\dots,D$), $X^i=\mathbb T^{M_i}, X=\prod\limits_{i=1}^D\mathbb T^{M_i}=\mathbb T^N$, в качестве меры на каждом торе берем стандартную меру Лебега, для $f(x)\in L^1(\mathbb T^{M_i})$ положим $ T^i_{n^i}f(x)=\sigma_{n^i}(f,x)=\int\limits_{\mathbb T^{M_i}}K_{n^i}(x^i-t^i)f(t^i)dt^i$ --- средние Фейера, и $K_{n^i}(\cdot)$ --- многомерное ядро Фейера с мультииндексом $n^i=(n^i_1,\dots,n^i_{M_i})$. В качестве направленного множества $A^i$ здесь выступает некоторое множество мультииндексов $\{(n^i_1,\dots,n^i_{M_i})\}$, возрастающих ограниченно (ограниченность понимается в смысле [1, гл.XVII, \S 3]), откуда следует, что соответствующий максимальный оператор $T^i$ имеет слабый тип $(1,1)$ (см. [1, лемма XVII.3.11]). Возрастание направленности мультииндексов $A^1\times\dots \times A^D$, где каждая направленность мультииндексов $A^i$ возрастает ограниченно, назовем $D$-ограниченным. Соответствующую сходимость результатов применения данных операторов назовем $D$-ограниченной. Известно, что средние Фейера ограниченно сходятся к исходной функции почти всюду если, например, исходная функция ограничена (см.[1, гл.XVII, \S 3]). Аналогичное верно и для некоторых других средних рядов Фурье, например, для средних Абеля-Пуассона (см. [1, гл.XVII, \S 3]), направленности для которых определяются аналогичным образом,  и средних Марцинкевича (см. [11] и [12]). 

Отметим, что в случае средних Абеля-Пуассона семейство операторов является не счетным, а континуальным, но в данном случае это не создает проблемы, так как можно переходить к счетному множеству индексов, плотному в исходном (например, рассматривая рациональные числа отрезка $[0,1]$), подобно тому, как это делается в [13, лемма 4.1.4]. $D$-кратными средними Марцинкевича функции $f\in L^1(\mathbb T^N)$ назовем средние вида 

$$
\sigma^{{\Phi}}_{\textbf n}( f,x)=\int\limits_{\mathbb T^N}\prod\limits_{j=1}^D K^{\Phi^j}_{ n^j}(x^j-t^j) f( t)d t
$$
где $K^{\Phi^j}_{ n^j}$ --- обычные ядра Марцинкевича, $t^j,x^j\in\mathbb T^{M_j}, t,x\in\mathbb T^N$. Подробно средние Марцинкевича изучены в [11]. Заметим, что при рассмотрении средних Марцинкевича каждая направленность $A^i$ состоит из мультииндексов вида $\{(n^i,\dots,n^i)\}\subset \mathbb N^{M_i}$, то есть состоит из семейства ``кубических'' мультииндексов.

Из вышесказанного и из теоремы 3 получаем следующий результат:

\begin{theorem} {Теорема 4.}  Пусть $ f \in L(\ln^{+}L)^{D-1} (\mathbb T^N)$. Тогда средние Фейера, Абеля-Пуассона и кратные средние Марцинкевича функции $f$ сходятся $D$-ограниченно к $f$ почти всюду на $\mathbb T^N$.

\end{theorem}

Другим примером служит обобщение теорем Лебега и Йессена-Марцинкевича-Зигмунда о дифференцировании в $\mathbb R^N$ неопределенных интегралов суммируемых функций. Подробнее об этом см. [14], там же приведены и дальнейшие приложения этого результата. Здесь же приведем только формулировку результата:

\begin{theorem} {Теорема 5.} Пусть $f \in L(\ln^{+}L)^{D-1} (\mathbb R^N) $. Тогда производная неопределенного интеграла функции $f$ для любой $D$-регулярной системы брусов существует и совпадает с $f$ почти всюду относительно меры Лебега.

\end{theorem}

\end{document}